\documentclass[11pt,twoside,english]{myamsart}
\advance\oddsidemargin by -1.0cm
\advance\evensidemargin by -1.0cm
\advance\topmargin by -1.0cm 
\usepackage{amssymb}
\usepackage{babel}
\usepackage{amstext}
\usepackage{amscd}   
\usepackage{epsfig}  
\usepackage{rotating}
\usepackage{graphicx}
\usepackage{psfrag}     
\theoremstyle{plain}
\newtheorem{cor}{Corollary}[section]
\newtheorem{lem}{Lemma}[section]
\newtheorem{thm}{Theorem}[section]
\newtheorem{prop}{Proposition}[section]
\theoremstyle{definition}
\newtheorem{exa}{Example}[section]
\newtheorem{NB}{Remark}[section]

\renewcommand{\labelenumi}{$(\theenumi)$}
\newcommand{\bdm}{\begin{displaymath}}
\newcommand{\edm}{\end{displaymath}}
\newcommand{\be}{\begin{equation}}
\newcommand{\ee}{\end{equation}}
\newcommand{\ba}[1]{\begin{array}{#1}}
\newcommand{\ea}{\end{array}}
\newcommand{\bea}[1][]{\begin{eqnarray#1}}
\newcommand{\eea}[1][]{\end{eqnarray#1}}
\newcommand{\btab}{\begin{tabular}}
\newcommand{\etab}{\end{tabular}}

\newcommand{\op}{\oplus}
\newcommand{\ox}{\otimes}

\newcommand{\ra}{\rightarrow}

\newcommand{\del}{\partial}

\newcommand{\C}{\ensuremath{\mathbb{C}}}
\newcommand{\R}{\ensuremath{\mathbb{R}}}

\newcommand{\Z}{\ensuremath{\mathbb{Z}}}

\newcommand{\vphi}{\ensuremath{\varphi}}    

\newcommand{\diag}{\ensuremath{\mathrm{diag}}}

\newcommand{\grad}{\ensuremath{\mathrm{grad}}}
\newcommand{\arctanh}{\ensuremath{\mathrm{arctanh}}}
\newcommand{\SO}{\ensuremath{\mathrm{SO}}}

\newcommand{\Orth}{\ensuremath{\mathrm{O}}}

\newcommand{\nms}{\!\!}%
\begin{document}
\def\haken{\mathbin{\hbox to 6pt{%
                 \vrule height0.4pt width5pt depth0pt
                 \kern-.4pt
                 \vrule height6pt width0.4pt depth0pt\hss}}}
    \let \hook\intprod
\setcounter{equation}{0}
%
%------ draw title page -----
%
\thispagestyle{empty}
%
%\hbox to \hsize{%
%  \vtop{} \hfill
%  \vtop{\hbox{PRELIMINARY VERSION}}}
%------------------------------
\date{\today}
%----------------------------------------------------------
\title[The geodesics of metric connections with vectorial torsion]{The 
geodesics of metric connections with vectorial torsion}
%----------------------------------------------------------
%
% author and address
%
%-------------------------------------------
%
\author{Ilka Agricola}
\author{Christian Thier}
%-------------------------------------------
\address{\hspace{-5mm} 
{\normalfont\ttfamily agricola@mathematik.hu-berlin.de}\newline
{\normalfont\ttfamily cthier@mathematik.hu-berlin.de}\newline
Institut f\"ur Mathematik \newline
Humboldt-Universit\"at zu Berlin\newline
Sitz: J.-von-Neumann-Haus Adlershof\newline
D-10099 Berlin, Germany}
%
%-----------------------------------------------------------
\thanks{Supported by the SFB 288 ``Differential geometry
and quantum physics'' of the DFG and the Junior Research Group
``Special geometries in mathematical physics'' of the Volkswagen 
Foundation.}  
%-----------------------------------------------------------
\subjclass[2000]{Primary 53 C 25; Secondary 81 T 30}
%-----------------------------------------------------------
%-----------------------------------------------------------
\begin{abstract}
%---------------
The present note deals with the dynamics of metric connections
with vectorial torsion, as already described by E.~Cartan in 1925.
We show that the geodesics of metric connections with  vectorial torsion 
defined by gradient vector fields coincide with the Levi-Civita geodesics 
of a conformally equivalent metric. By pullback, this yields a systematic 
way of constructing invariants of motion for such connections from 
isometries of the conformally equivalent metric, and we explain in as much
this result generalizes the Mercator projection which maps sphere loxodromes
to straight lines in the plane. An example shows that Beltrami's
theorem fails for this class of connections. We then study the system of 
differential
equations describing geodesics in the plane for vector fields which are not 
gradients, and show among others that the Hopf-Rinow theorem does
also not hold in general. 
\end{abstract}
%-------------
\maketitle
%----------------
%\tableofcontents
%----------------
\pagestyle{headings}
%
%
%-------------- body of the document ------------------------------------------
%
%---------------------------------------------------------------------------- 
\section{Introduction}\noindent
%----------------------------------------------------------------------------
%
The present note deals with the dynamics of metric connections
on a given Riemannian manifold $(M^n,g)$ of the form
\bdm
\nabla_X Y\ =\ \nabla^g_XY+g(X,Y)V-g(V,Y)X,
\edm
where $V$ denotes a fixed vector field on $M$ and $\nabla^g$ is
the usual Levi-Civita connection. This is one of the three basic
types of metric connections, as already described by E.~Cartan in 
\cite{Cartan25a}, and will be called by us a metric connection 
\emph{with vectorial torsion}. The case of a surface ($n=2$) is special in 
as much that any metric connection has to be of this type. In fact,
classical topics of surface theory like the Mercator projection which
maps loxodromes on the sphere  to straight lines in the plane can be 
understood in a different light
with their help. We will show that the geodesics of metric connections with  
vectorial torsion defined by gradient vector fields coincide with the 
Levi-Civita geodesics of a conformally equivalent metric. By pullback, 
this yields a systematic way of constructing invariants of motion for such 
connections from  isometries of the conformally equivalent metric. 
Furthermore, the discussion of the 
catenoid will show that there exist geodesic mappings for this larger
class of connections between surfaces of constant and non-constant
Gaussian curvature, hence Beltrami's theorem becomes wrong in this
situation. In the last Section, we will treat the euclidian plane
with an arbitrary vector field in great detail and show
that the behaviour of geodesics is, already in this low dimension,
dominated by a highly non-trivial system of ordinary differential
equations. In particular, the Hopf-Rinow theorem can fail, too.

\smallskip\noindent
For some reasons, these connections have not attracted as much
attention in the past  as we believe they deserve. Correspondingly,
an overview over the existing literature (that we are aware of) is 
quickly given. In \cite{Tricerri&V1}, Tricerri and Vanhecke were led
to the study of such connections in the context of the classification
problem of homogeneous structures on manifolds. They showed
that if $M$ is connected, complete, and simply connected and
$V$ is parallel, i.\,e.~$\nabla V=0$, then $(M,g)$ has to be
isometric to hyperbolic space.
Vicente Miquel studied in \cite{Miquel82} and \cite{Miquel01} the
growth of geodesic balls of such connections, but did not
investigate the detailed shape of  geodesics. Connections with vectorial
torsion on spin manifolds may also play a role in superstring theory (see 
\cite{Agricola&F03a} and the literature cited therein), but this aspect
will not be discussed in the present paper.
%
%\smallskip\noindent
Both authors wish to thank Thomas Friedrich and Pawe{\l} Nurowski for 
valuable discussions and suggestions.

\smallskip\noindent
We finish this introduction with  a short review of the possible classes of 
metric connections. For details, we refer to \cite{Cartan25a}, 
\cite{Tricerri&V1} and \cite{Agricola&F03a}.
In any  point $p$ of an $n$-dimensional Riemannian manifold $(M^n, g)$,
the difference between its Levi-Civita connection $\nabla^g$ and any linear 
connection $\nabla$ is  a $(2,1)$-tensor $A$, 
\bdm 
\nabla_X Y\ =\ \nabla^g_X Y + A(X,Y),\quad X,Y \in T_p M\,.
\edm
We identify $T_p M$ with $\R^n$, on which the real orthogonal group  
$\Orth(n,\R)$ then acts in the standard way.
The connection $\nabla$ is metric if and only if the difference tensor $A$
(viewed as a $(3,0)$-tensor)  belongs in every point to the 
$n^2(n-1)/2$-dimensional space
\bdm
\mathcal{A}^g\ :=\ \R^n\ox\wedge^2 \R^n \ = \ \{A \in\ox^3 \R^n \ | 
\ A(X,V,W) +  A(X,W,V) \ =\ 0\} \, .
\edm
\begin{prop}\label{classes}
%---------------------------
For $n\geq 3$, the space $\mathcal{A}^g$ of possible metric 
difference tensors splits under $\Orth(n,\R)$ into the sum of three 
irreducible representations,
$\mathcal{A}^g \cong \R^n\op \wedge^3 \R^n \op \mathcal{A}'$.
For $n=2$, $\mathcal{A}^g$ is an irreducible $\Orth(2,\R)$-module
isomorphic to $\R^2$.
\end{prop}
\noindent
The torsion of $\nabla$ is then constructed from $A$ by the
usual formula $T(X,Y)=A(X,Y)-A(Y,X)$.
It has become customary to call metric connections such that
$A\in \wedge^3 \R^n$  \emph{metric connections with totally skew-symmetric
torsion}, and we shall use the denomination \emph{metric connections with 
vectorial torsion} in case $A\in\R^n$. The correspondence can be made
explicit: for a vector field $V$ on $M$, $A$ and $T$ are then given by 
\bdm
A(X,Y)\ =\ g(X,Y)V-g(V,Y)X,\quad
T(X,Y,Z)\ =\ g\big(g( V, X)Y- g(V, Y)X,Z\big).
\edm
The metric connection $\nabla$ has the same geodesics as the Levi-Civita
connection precisely if it has totally skew-symmetric torsion. 
%
%---------------------------------------------------------------------------
\section{Basic remarks on geodesics of metric connections with vectorial 
torsion}\label{general-remarks}\noindent
%---------------------------------------------------------------------------
%
Let $(M,g)$ be a Riemannian manifold, $V$ a vector field on $M$ and
$\nabla$ the metric connection with vectorial torsion defined by $V$.
A curve  $\gamma(t)$  is a geodesic of $\nabla$ if it satisfies the
differential equation 
\bdm
\nabla^g_{\dot{\gamma}}\dot{\gamma} + g(\dot{\gamma},\dot{\gamma})V -
g(V,\dot{\gamma})\dot{\gamma}\ = \ 0\,.
\edm
Taking the scalar product of this equation with $\dot{\gamma}$ yields
$g(\nabla^g_{\dot{\gamma}}\dot{\gamma}, \dot{\gamma})=0$, that is, 
$\dot{\gamma}$
has constant length $E>0$, which reflects of course just the fact that
$\nabla$ was metric. Hence, the geodesic equation can be written
\be\label{geod-eq}
\nabla^g_{\dot{\gamma}}\dot{\gamma} + E^2\,V -
g(V,\dot{\gamma})\dot{\gamma}\ = \ 0\,.
\ee
In fact, there are qualitatively two cases to be distinguished. 
If $\dot{\gamma}$ is parallel to $V$ at the origin,  
$\dot{\gamma}(0)=\alpha\cdot V(\gamma(0))$, we conclude that 
$\nabla^g_{\dot{\gamma}(0)}\dot{\gamma}(0)=0$ and $\gamma(t)$ coincides
locally with a classical geodesic of the Levi-Civita connection.
In particular, a $\nabla$-geodesic which stays parallel to $V$ for
all times is exactly a $\nabla^g$-geodesic. Generic $\nabla$-geodesics
are those for which $\dot{\gamma}$ is never parallel to $V$; their shape
will be qualitatively very different from that of their Levi-Civita 
``cousins''. First, we express the  curvature of a geodesic curve
through $V$ and $\dot{\gamma}$:
\begin{lem}\label{gen-curv}
%--------------------------
The Riemannian geodesic curvature $\kappa$ of a $\nabla$-geodesic $\gamma$ 
is given by
$\kappa^2=||V||^2 - g(V,\dot{\gamma})^2/E^2$ for general $V$,
and by $\kappa^2=- E^{-2}\,d/dt\, g(V,\dot{\gamma})$ if $V$ is a Killing vector
field.
\end{lem}
\begin{proof}
Observe that $\ddot{\gamma}=\nabla^g_{\dot{\gamma}}\dot{\gamma}$, and
$\kappa=||\ddot{\gamma}||/E^2$ in our normalization. Differentiating the
geodesic equation (\ref{geod-eq}) once more and taking the scalar product
with $\dot{\gamma}$,  we obtain
\bdm
g(\dddot{\gamma},\dot{\gamma})\ =\ E^2 g(V,\ddot{\gamma})\ =\ 
E^2[g(V,\dot{\gamma})^2-E^2||V||^2 ]\,.
\edm
But from $g(\dot{\gamma},\ddot{\gamma})=0$, we conclude that
$g(\ddot{\gamma},\ddot{\gamma})+g(\dot{\gamma},\dddot{\gamma})=0$, hence the
first claim follows. Suppose now that $V$ is a Killing vector field,
$g(\nabla^g_X V,Y)+g(\nabla^g_Y V,X)=0$ for all $X$ and $Y$.
This implies in particular 
$g(\dot{\gamma},\dot{V})=g(\dot{\gamma},\nabla^g_{\dot{\gamma}}V)=0$,
and so
\bdm
\frac{d}{dt} g(V,\dot{\gamma})\ =\ g(V,\ddot{\gamma})\ =\
g(V,\dot{\gamma})^2-E^2||V||^2\ =\ -E^2\kappa^2. \qedhere
\edm
\end{proof}
\noindent
In examples, this allows to determine a differential equation for $\kappa$,
which one then needs to study in detail. Even without this, the lemma
is useful for proving results of the following type: periodic
$\nabla$-geodesics of Killing vector fields have to be of non-generic type.
More precisely:
\begin{cor}
%----------
For a Killing vector field $V$, any periodic $\nabla$-geodesic is
automatically a Levi-Civita geodesic and, up to a constant, an integral curve
of $V$.
\end{cor}
\begin{proof}
%-------------
If the curve $\gamma$ is periodic, the functions $g(V,\dot{\gamma})$ and
$d/dt \,g(V,\dot{\gamma})$ are periodic too. But since the latter
has to be non-positive by Lemma~\ref{gen-curv}, we conclude that
$d/dt \,g(V,\dot{\gamma})=0$, hence $\nabla^g_{\dot{\gamma}}\dot{\gamma}=0$.
Furthermore, the geodesic equation  yields $E^2V=c\,\dot{\gamma}$,
where $c$ is the constant $g(V,\dot{\gamma})$, and $(c/E^2)\gamma$ becomes
an integral curve of $V$. 
\end{proof}
\noindent
In the next step, we show that isometries leaving $V$ invariant generate
symmetries of the $\nabla$-geodesics.
\begin{prop}\label{isometries}
%-----------------------------
If $X$ is a Killing vector field commuting with $V$, its flow
$\Phi_s$ maps $\nabla$-geodesics to $\nabla$-geodesics.
\end{prop}
\begin{proof}
%------------
Let $\gamma(t)$ be a $\nabla$-geodesic, and set 
$\gamma^*(t):=\Phi_s\big(\gamma(t)\big)$, for which one has
$\dot{\gamma}^*=d\Phi_s(\dot{\gamma})$. Since $\Phi_s$ consists of isometries, 
$||\dot{\gamma}||=||\dot{\gamma}^*||=E$, and the invariance of $V$ under
$\Phi_s$ implies $g(V,\dot{\gamma})=g(V,\dot{\gamma}^*)$.
In addition, $d\Phi_s(\nabla^g_{\dot{\gamma}}\dot{\gamma})=
\nabla^g_{d\Phi_s(\dot{\gamma})} d\Phi_s(\dot{\gamma})$, hence applying 
$d\Phi_s$ to the geodesic equation for $\gamma(t)$ yields that
$\gamma^*(t)$ satisfies the geodesic equation, too. Notice that by the same
argument, the first curvatures of $\gamma$ and $\gamma^*$ coincide. 
\end{proof}
%
%\noindent
%In Section \ref{gradVFs}, we will construct  a  geodesic
%mapping for metric connections with vectorial torsion that is not an
%isometry leaving $V$ invariant. 
%
\begin{NB}\label{no-noether}
%---------------------------
In Section \ref{gradVFs}, we will construct  a  geodesic
mapping for metric connections with vectorial torsion that is not an
isometry leaving $V$ invariant. 
Also, the main problem in the study of the geodesic equation is that 
Proposition 
\ref{isometries} does \emph{not} yield invariants of motion of Noether type:
if $\gamma$ is a $\nabla$-geodesic and $W$ a vector field  whose flow
$\Phi_s$ consists of isometries and leaves $V$ invariant,
the quantity $g(W,\dot{\gamma})$ is \emph{not} a first integral
(easy examples will again be discussed in  Section \ref{gradVFs}).
This reflects the fact that  the geodesic equation is not the
Euler-Lagrange equation of a Lagrange function 
$\mathcal{L}(\dot{\gamma})=g_{ij}\dot{x}^i\dot{x}^j/2- V(\gamma)$
for some potential function $V$.
\end{NB}

%---------------------------------------------------------------------------
\section{Geodesics for gradient vector fields and conformal mappings}
\label{gradVFs}
\begin{figure}
\bdm
\includegraphics[width=8cm, bbllx=50, bblly=230, bburx=560,bbury=590]{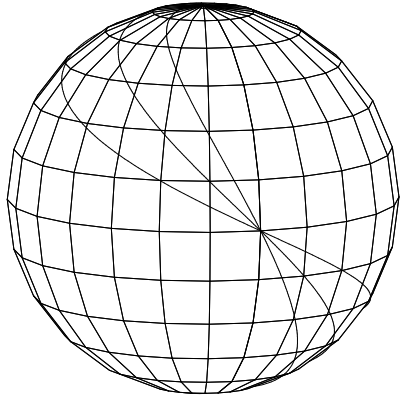}
\edm
\caption{Loxodromes on the sphere.}\label{sphaere}
\end{figure}

\begin{figure}
\begin{center}
\begin{minipage}[b]{.46\linewidth}
\epsfig{figure=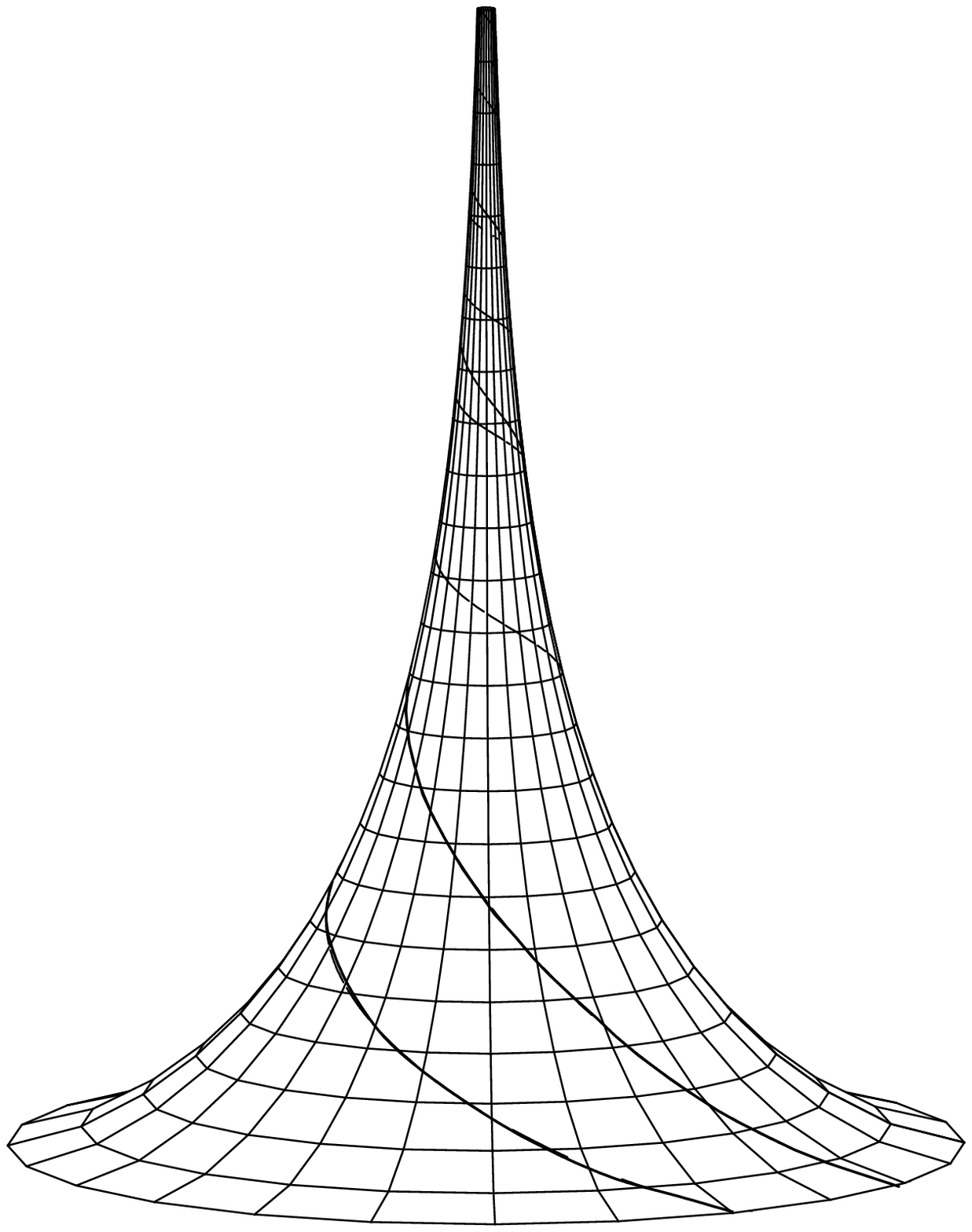,width=\linewidth}
\end{minipage}\hfill
%---------------------------------
\begin{minipage}[b]{.46\linewidth}
%\centering
\epsfig{figure=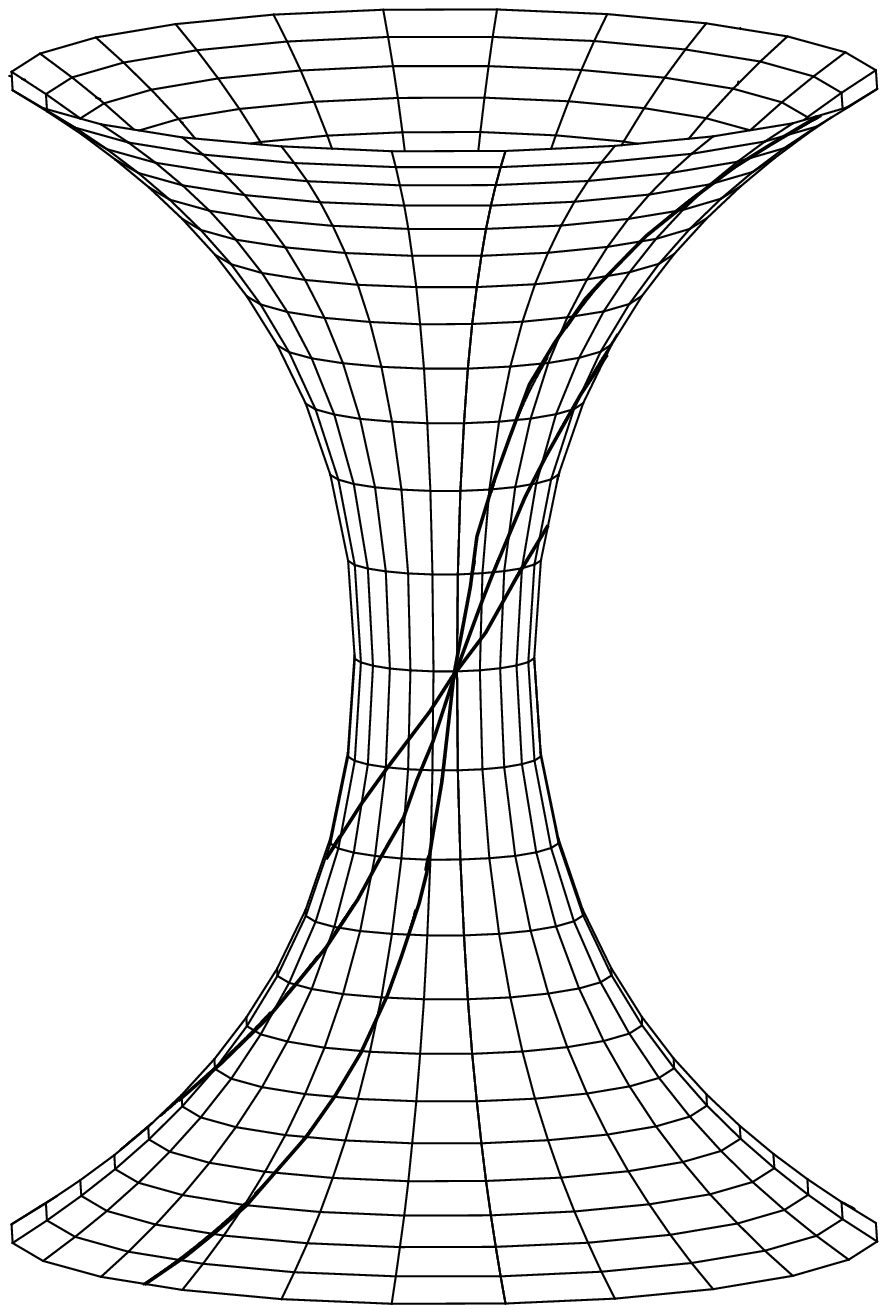, width=\linewidth}
\end{minipage}
\end{center}
\caption{Loxodromes on the pseudosphere and the  catenoid.}
\label{pseudo&kat}
\end{figure}

\noindent
%---------------------------------------------------------------------------
%
The motivation for this section was  Cartan's example of a metric connection
with vectorial torsion on the sphere. Up to our knowledge, this is the only
instance where some geodesics of metric connections with vectorial torsion
are described in the literature.
\begin{exa}[Cartan's example]\nocite{Cartan24a}
%----------------------------------------------
In \cite[\S\,67, p.\,408--409]{Cartan23a}, Cartan describes the two-dimensional
sphere with its flat metric connection, and observes (without proof) that 
``on this manifold, the straight lines are the \emph{loxodromes}, which 
intersect the
meridians at a constant angle. The only straight lines realizing shortest
paths are those which are normal to the torsion in every point: \emph{these
are the meridians}\footnote{Sur cette vari\'et\'e, les lignes droites sont
les \emph{loxodromies}, qui font un angle constant avec les m\'eridiennes.
Les seules lignes droites qui r\'ealisent les plus courts chemins sont celles
qui sont normales  en chaque point \`a la torsion : \emph{ce sont les 
m\'eridiennes.} loc.\,cit.}''. 
\end{exa}
\noindent
This suggests that  there exists a class of  metric connections on surfaces 
of revolution whose geodesics admit a generalization of Clairaut's theorem,
yielding loxodromes in the case of the flat connection. Furthermore,
it is well known that the Mercator projection maps loxodromes to 
straight lines in the plane (i.\,e., Levi-Civita geodesics of the euclidian
metric), and that
this mapping is conformal. Theorem~\ref{thm-grad-VF} provides the right
setting to understanding both effects, as explained in the introduction.
First, let us recall without proof a  standard formula for a conformal change 
$\tilde{g}=e^{2\sigma}g$ of a metric $g$ that we will need in the sequel: 
\begin{prop}\label{conf-change}
%-----------
The Levi-Civita connection behaves as follows under conformal change of 
the metric:
\bdm
\nabla^{\tilde{g}}_X Y\,=\, \nabla^g_X Y + X(\sigma)Y+Y(\sigma)X
-g(X,Y)\grad (\sigma).
\edm
\end{prop}
\begin{thm}\label{thm-grad-VF}
%-----------------------------
Let $\sigma$ be a function on the Riemannian manifold $(M,g)$, 
$\nabla$ the metric connection with vectorial torsion defined by
$V=-\grad(\sigma)$, and consider the conformally equivalent metric
$\tilde{g}=e^{2\sigma}g$. Then:
\begin{enumerate}
\item Any $\nabla$-geodesic $\gamma(t)$ is, up to
a reparametrisation $\tau$, a $\nabla^{\tilde{g}}$-geodesic, and 
the function $\tau$ is the unique solution of the differential
equation $\ddot{\tau}+\dot{\tau}\dot{\sigma}=0$, where we set
$\sigma(t):=\sigma\circ\gamma\circ\tau (t)$;
\item If $X$ is a Killing field for the metric $\tilde{g}$, the function
$e^{\sigma}g(\dot{\gamma}, X)$ is a constant of motion for the
$\nabla$-geodesic $\gamma(t)$.
\end{enumerate}
\end{thm}
\begin{proof}
Assume that $\gamma: I \ra M$ is a $\nabla$-geodesic, i.\,e.~it satisfies
\bdm
\nabla^g_{\dot{\gamma}}\dot{\gamma}\, =\, - E^2\,V +
g(V,\dot{\gamma})\dot{\gamma}\, =\, E^2\,\grad(\sigma) -
g(\grad(\sigma),\dot{\gamma})\dot{\gamma}.
\edm
Then Proposition~\ref{conf-change} implies
\be\label{conf-deriv}
\nabla^{\tilde{g}}_{\dot{\gamma}(t)}\dot{\gamma}(t) \,=\, 
2\dot{\gamma}(\sigma)\dot{\gamma} - g(\grad(\sigma),\dot{\gamma})\dot{\gamma}
\, =\, +\, g(\grad(\sigma),\dot{\gamma})\dot{\gamma}\, =\,
\frac{d}{dt} (\sigma\circ\gamma)(t)\cdot \dot{\gamma}(t).
\ee
We claim that there exists a reparametrisation 
$\gamma^*(t):=\gamma\circ\tau(t)$ of the curve $\gamma(t)$
which satisfies $\nabla^{\tilde{g}}_{\dot{\gamma}^*}\dot{\gamma}^*=0$.
For a still arbitrary function $\tau:I\ra I$, we have
$\dot{\gamma}^*(t)\ =\ \dot{\tau}(t)\cdot\dot{\gamma}(\tau(t))$,
hence
\bea[*]
\nabla^{\tilde{g}}_{\dot{\gamma}^*(t)}\dot{\gamma}^*(t) & = &
\nabla^{\tilde{g}}_{\dot{\gamma}^*(t)}\big(\dot{\tau}(t)\cdot\dot{\gamma}
(\tau(t))\big)\ =\
\dot{\gamma}^*(\dot{\tau})\cdot\dot{\gamma}(\tau(t))+\dot{\tau}(t)\,
\nabla^{\tilde{g}}_{\dot{\gamma}^*(t)}\dot{\gamma}(\tau(t))\\ 
& = & \dot{\tau}\dot{\gamma}(\dot{\tau})\cdot\dot{\gamma} +\dot{\tau}^2 \,
\nabla^{\tilde{g}}_{\dot{\gamma}(\tau(t))}\dot{\gamma}(\tau(t))
\ =\  \left[\dot{\tau}\dot{\gamma}(\dot{\tau}) +\dot{\tau}^2 \,
\frac{d}{ds}\bigg|_{s=\tau(t)} \nms\nms\nms\nms\nms\nms(\sigma\circ\gamma)(s)
 \right]\dot{\gamma}(\tau(t)),
\eea[*]
where we used equation~(\ref{conf-deriv}) in the last step. It remains to be 
shown that $\tau$ can be chosen such that the expression in parentheses 
vanishes.  For this, we first rewrite this expression as
\bdm
\dot{\tau}(t)\frac{d}{ds}\bigg|_{s=\tau(t)} \nms\nms\nms\nms\nms\nms
(\dot{\tau}\circ\gamma)(s) +\dot{\tau}^2(t) \,
\frac{d}{ds}\bigg|_{s=\tau(t)} \nms\nms\nms\nms\nms\nms\nms
(\sigma\circ\gamma)(s) \ =\
\frac{d}{ds}\bigg|_{s=t} (\dot{\tau}\circ\gamma\circ\tau)(s)
 +\dot{\tau} \,
\frac{d}{ds}\bigg|_{s=t} \nms\nms\nms(\sigma\circ\gamma\circ\tau)(s),
\edm
leading to the differential equation 
\be\label{reparam}
\ddot{\tau}(t)+\dot{\tau}(t)\dot{\sigma}(t)\ =\ 0
\ee
for the function $\tau$. Here, we defined 
$\sigma(t):=(\sigma\circ\gamma\circ\tau)(t)$ and viewed $\dot{\tau}$ as a
function on $\gamma(I)\subset M$ if necessary. It cannot be solved explicitely,
since $\sigma$ depends implicitly on $\tau$; but by Picard-Lindel\"of, 
equation~(\ref{reparam}) does  always admit a solution. Observe that,
formally, the differential equation can be integrated once, yielding
the relation
$%\bdm
\dot{\tau}(t)\ =\ e^{-\sigma(t)}\,.
$%\edm
%
%------------------------------------------------------------
\begin{figure}%\label{drehflaeche}
\bdm
\psfrag{x}{$x$}\psfrag{y}{$y$}\psfrag{z}{$z$}
\psfrag{a}{$h(s)$}\psfrag{r}{$r(s)$}
\psfrag{p}{$\varphi$}
\psfrag{v1}{$\nu_1$}\psfrag{v2}{$\nu_2$}
\includegraphics[width=6.6cm]{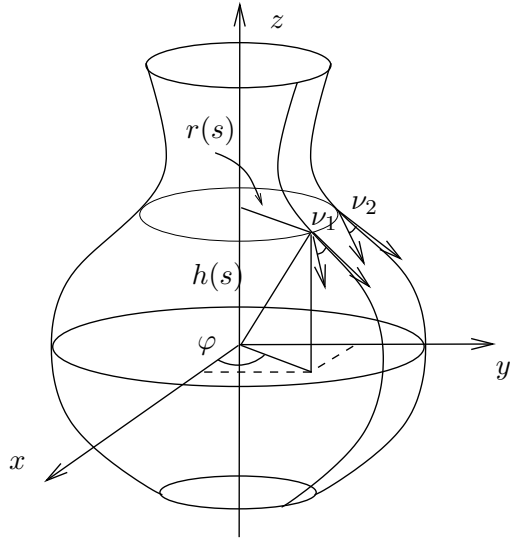}
\edm
\caption{Surface of revolution generated by a curve $\alpha$.}
\label{drehflaeche}    
\end{figure}
%-------------------------------------------------------------
%
%
For the second part of the proof, let us now assume that $X$ is a Killing 
vector field relatively to the metric $\tilde{g}=e^{2\sigma}g$.
Thus, the $\nabla^{\tilde{g}}$-geodesic $\gamma^*$ satisfies 
$\tilde{g}(\dot{\gamma}^*,X\circ\gamma^*)=\mathrm{const}=:c$.
By the integrated form of  equation~(\ref{reparam}), this implies
\bdm
c\ = \ e^{2\sigma(t)} g \big(\dot{\tau}(t) \dot{\gamma}(\tau(t)),
X\big|_{\gamma(\tau(t))}\big)\ =\ e^{2\sigma(t)} e^{-\sigma(t)}
g \big(\dot{\gamma}(\tau(t)), X\big)\ =\
e^{\sigma} g(\dot{\gamma},X). \qedhere
\edm
%
%This finishes the proof of Theorem \ref{thm-grad-VF}.
\end{proof}
\noindent
Since one is most often only interested in the set of points of a curve and 
not in the  parametrisation itself, it is rarely necessary
to determine $\tau$ explicitely.
\begin{exa}[Loxodromes and Mercator projection]
%----------------------------------------------
%
We discuss Cartan's example in the light of Theorem \ref{thm-grad-VF}.
Let $\alpha=(r(s),h(s))$ be a curve in natural parametrisation,
and $M(s,\vphi)=(r(s)\cos\vphi,r(s)\sin\vphi,h(s))$ the surface of
revolution generated by it. The first fundamental form is then
$g=\diag(1,r^2(s))$, and we fix the orthonormal frame
$e_1=\del_s$, $e_2=(1/r)\del_\vphi$ with dual $1$-forms
$\sigma^1=ds$, $\sigma^2=r\,d\vphi$. We convene to call two tangential vectors
$v_1$ and $v_2$ parallel if their angles $\nu_1$ and $\nu_2$ with
the meridian through that point coincide (see Figure \ref{drehflaeche}).
Hence $\nabla e_1=\nabla e_2=0$, and the connection $\nabla$ is flat.
But for a flat connection, the torsion $T$ is can be derived from
$d\sigma^i(e_j,e_k)=\sigma^i(T(e_j,e_k))$. Since $d\sigma^1=0$ and
$d\sigma^2=(r'/r)\sigma^1\wedge\sigma^2$, one obtains 
\bdm
T(e_1,e_2)\ =\ \frac{r'(s)}{r(s)}\,e_2 \ \text{ and }\
V\ =\ \frac{r'(s)}{r(s)}\,e_1\ =\ -\grad \big(- \ln r(s)\big).
\edm
Thus, the metric connection $\nabla$ with vectorial torsion $T$ is
determined by the gradient of the function $\sigma:=- \ln r(s)$. 
By Theorem \ref{thm-grad-VF}, we conclude that its geodesics are
the Levi-Civita geodesics of the conformally equivalent metric
$\tilde{g}=e^{2\sigma}g=\diag(1/r^2, 1)$. This coincides with the standard
euclidian metric if one performs the change of variables
$x=\vphi$, $y=\int ds/r(s)$. For example, the sphere is obtained
for $r(s)=\sin s, h(s)=\cos s$, hence $y=\int ds/\sin s = \ln\tan(s/2)$
 ($|s|<\pi/2$), and
this is precisely the coordinate change of the Mercator projection.
Furthermore, $X=\del_\vphi$ is a Killing vector field for $\tilde{g}$, 
hence the second part of Theorem \ref{thm-grad-VF}
yields for a $\nabla$-geodesic $\gamma$ the invariant of motion
\bdm
\mathrm{const}\ =\ e^\sigma g(\dot{\gamma},X)\ =\
\frac{1}{r(s)} g(\dot{\gamma}, \del_\vphi)\ = \ g(\dot{\gamma},e_2),
\edm
and this is just the cosine of the angle between $\gamma$ and a parallel
circle. This shows that $\gamma$ is a loxodrome on $M$, as claimed
(see Figure \ref{sphaere} for loxodromes on the sphere). In the same way,
one obtains a ``generalized Clairaut theorem'' for any gradient
vector field on a surface of revolution. Figure \ref{pseudo&kat}
illustrates loxodromes on a pseudosphere and  a catenoid. For the
pseudosphere, one chooses
\bdm
r(s)\ =\ e^{-s},\quad h(s)\ =\ \arctanh \sqrt{1-e^{-2s}} - \sqrt{1-e^{-2s}},
\edm
hence $V=-e_1$ and $\nabla V=0$, in accordance with the results by
\cite{Tricerri&V1} cited in the introduction. Notice that $X$ is also
a Killing vector field for the metric $g$ and does commute with $V$;
nevertheless, $g(\dot{\gamma},X)$ is \emph{not} an invariant of motion
(compare this to Proposition~\ref{isometries} and Remark~\ref{no-noether}).
\end{exa}
%
%--------------------------------------------------------------
%--------------------------------------------------------------
\begin{NB}
%---------
The catenoid is interesting for another reason:
since it is a minimal surface, the Gauss map to the sphere is a 
conformal mapping, hence it maps loxodromes to loxodromes. Thus,
Beltrami's theorem (``If a portion of a surface $S$ can be mapped 
LC-geodesically onto a portion of a surface $S^*$ of constant Gaussian
curvature, the Gaussian curvature of $S$ must also be constant'',
see for example \cite[\S 95]{Kreyszig2}) does \emph{not} hold for metric
connections with vectorial torsion -- the sphere is a a surface
of constant Gaussian curvature, but the catenoid is not.
\end{NB}
%
%------------------------------------ 
%--------------------------------

\begin{NB}
%---------
The unique flat metric connection $\nabla$ does not have to be of vectorial 
type.
For example, on the compact Lie group $\SO(3,\R)$, its torsion is 
a $3$-form: Fix an orthonormal basis $e_1,e_2,e_3$ with
commutator relations $[e_1,e_2]=e_3, [e_2,e_3]=e_1$ and $[e_1,e_3]=-e_2$.
Cartan's structural equations then read $d\sigma^1=\sigma^2\wedge\sigma^3,
 d\sigma^2=-\sigma^1\wedge\sigma^3,d\sigma^3=\sigma^1\wedge\sigma^2$,
from which we deduce $T=2A=\sigma^1\wedge\sigma^2\wedge\sigma^3$.
In particular, $\nabla$ has the same geodesics as $\nabla^g$.
\end{NB}
\begin{NB}
%---------
The spirit of Theorem \ref{thm-grad-VF} reflects a well-known fact from Weyl 
geometry. Recall that a \emph{Weyl manifold}  is a manifold $M^n$ 
with a conformal
structure $c$ and a torsion-free connection on the $\mathrm{CO}(n)$-reduction
of the frame bundle of $M^n$. For fixed $c$, the set of all such connections 
is an affine space modelled over the vector space of all $1$-forms. 
The choice of a metric $g$ in the conformal class $c$ simply identifies the 
connection  with a $1$-form $\theta\in\Omega^1(M^n)$. For a second 
metric $g'=e^{-2f}g$ in the conformal class $c$, the 
characteristic $1$-form $\theta'$ is then given by $\theta'=\theta+df$.
\end{NB}
%
%---------------------------------------------------------------------------- 
\section{Geodesics of vectorial connections in the plane}
\label{geod-plane}\noindent
%---------------------------------------------------------------------------- 
%
We consider the plane $M=\R^2$ endowed with its standard euclidian
metric $g=dx^2+dy^2$, and the metric connection $\nabla$ with vectorial
torsion defined by the vector field $V=f(x,y)\,\del_x+ g(x,y)\,\del_y$, where
$f$ and $g$  may be any smooth functions of $x$ and $y$. Then one computes 
for the standard orthonormal basis $e_1=\del_x,\,e_2=\del_y$ that 
\bdm
\nabla_{e_1}e_1\,=\, g\,e_2,\quad\nabla_{e_1}e_2\,=\,-g\,e_1,\quad
\nabla_{e_2}e_1\,=\,-f\,e_2,\quad\nabla_{e_2}e_2\,=\,f\,e_1,
\edm
and hence the connection form $\omega_{12}$ and the curvature form
$d\omega_{12}$ are given by
\bdm
\omega_{12}\ =\ g(x,y)\,dx-f(x,y)\,dy,\quad
d\omega_{12}\ =\ -(\del_xf+\del_y g)\,dx\wedge dy.
\edm
The geodesics $\gamma(t)=(x(t),y(t))$ of $\nabla$ are the 
solutions of the differential equations
\bdm
\ddot{x}-g(x,y)\,\dot{x}\dot{y} +f(x,y)\,\dot{y}^2 \,=\,0,\quad
\ddot{y}-f(x,y)\,\dot{x}\dot{y}+ g(x,y)\,\dot{x}^2\,=\,0.
\edm
Since a plane curve admits a notion of curvature endowed with a sign
and is, up to euclidian motions, uniquely determined by it, the general
considerations of Section~\ref{general-remarks} can be 
strengthened in this case. The curvature $\kappa$ is defined through
the requirement $\ddot{\gamma}= \kappa \cdot i\cdot \dot{\gamma}$,
hence we obtain at once: 
\begin{lem}\label{geod-plane-1}
%--------------------------
The equations for a geodesic $\gamma$ of the metric connection $\nabla$ with 
vectorial torsion defined by $V=f(x,y)\,\del_x+ g(x,y)\,\del_y$ are
\bdm
\ddot{x}\ =\ -\kappa\,\dot{y},\quad \ddot{y}\ =\ \kappa\,\dot{x},
\edm
and $\kappa=f(x,y)\,\dot{y}- g(x,y)\,\dot{x}$ is the curvature of the plane
curve  $\gamma$.
\end{lem}
\noindent
If $\del_y f=\del_x g$, the vector field $V$ is (locally) a 
gradient vector field and its geodesics are best treated with the methods of
Section~\ref{gradVFs}. The form of the differential equations does 
already suggest that the system will, in general, not admit an easy solution. 
In fact, even a discussion of its qualitative behaviour turns out to be rather
difficult.

Yet, there is a second case which turns out to be manageable,
namely, the flat case. By the general form of $d\omega_{12}$ given above,
the connection $\nabla$ will be flat if there exists a smooth 
function $p(x,y)$ such that
$f=\del_y p$ and $g=-\del_x p$. The time derivative of $p$ is then
$\dot{p}=\del_xp\,\dot{x}+\del_yp\,\dot{y}=-g\,\dot{x}+f\,\dot{x}=\kappa$,
so the two geodesic equations of Lemma \ref{geod-plane-1} can be
restated in one single equation for the complex variable $z=x+iy$,
\bdm
\ddot{z}\ =\ i\,\dot{p}\dot{z}.
\edm
If $\dot{z}\neq 0$, we divide both sides by $\dot{z}$ and 
obtain a second invariant of motion (besides $E$)
\bdm
\frac{d}{dt}(\ln \dot{z}-ip)\ =\ \mathrm{const}\in \C.
\edm
By choosing the constants appropriatly, we finally arrive at
\bdm
\dot{z}\ =\ z_0 e^{ip(x,y)}\quad \text{for some } z_0\in\C.
\edm
For given $p$ and $z_0$, the solution  of this differential equation
may now be drawn with any standard ODE computer package. 
Notice that, strictly speaking, only the modulus of this equation
yields a new invariant of motion; taking its absolute value,
we arrive at $ |\dot{z}|=|z_0|=1$ if we assume---as we always did---that
$\gamma$ is given in natural parametrisation.  Obviously, this
is not an invariant of Noether type.
We give
two examples which illustrate how the behaviour of geodesics can differ 
for different functions $p$.
\begin{exa}
%----------
For $p=-(x^2+y^2)/2$, one has $f=-y,g=+x$, and $V$ is the winding vector 
field in the plane. The
geodesic  equations  can be written 
\bdm%\tag{$**$}
\ddot{x}\,=\,k\,\dot{y},\quad \ddot{y}\,=\,-k\,\dot{x},
\edm
where $k=-\kappa =x\dot{x}+y\dot{y}=d/dt\, ||\gamma(t)||^2/2$:
the curvature of $\gamma$ coincides, up to sign, with
half the derivative of the distance between the point $\gamma(t)$
and the origin. This illustrates that the origin plays indeed a 
special role for this vector field.  $V$ is the Killing
vector field generating rotations $r_{\varphi}$ by any angle $\varphi$ 
around the origin, hence Proposition~\ref{isometries} implies that
the curve $r_{\varphi}(\gamma)$ is again a geodesic if $\gamma$ was one.
%
%\begin{lem}
%%----------
%For $E=1$, the function $k$ satisfies the differential equations
%%
%\bdm
%k\ddot{k}\ =\ \dot{k}(\dot{k}-1)-k^4,\quad 
%\dot{k}(\dddot{k}+4k^2\dot{k})\ =\ \ddot{k}(\ddot{k}+k^3)
%\edm 
%%
%and the initial conditions
%\bdm
%k_0\,=\,x_0\dot{x}_0+y_0\dot{y}_0,\quad
%\dot{k}_0\,=\, 1+ x_0\ddot{x}_0+y_0\ddot{y}_0,\quad
%\ddot{k}_0\,=\, x_0\dddot{x}_0+y_0\dddot{y}_0.
%\edm
%%
%\end{lem}
%%
%\begin{proof}
%%------------
%Explicit differentiation of $k$ yields
%$\dot{k}=1+ x\ddot{x}+y\ddot{y}$ and $\ddot{k}= x\dddot{x}+y\dddot{y}$.
%On the other hand, we deduce from the geodesic equations that
%$\dddot{x}=\dot{k}\ddot{x}/k-k^2\dot{x}$ and
%$\dddot{y}=\dot{k}\ddot{y}/k-k^2\dot{y}$. Substituting these results into
%the expression for $\ddot{k}$ yields the first differential equation.
%The second differential equation is then obtained by a routine
%computation.
%\end{proof}
%%
%\noindent
%In particular, for a geodesic through $(0,0)$, one has
%$k_0=0, \dot{k}_0=1$ and $\ddot{k}_0=0$, hence the second order equation for 
%$k$ has a singularity at $t=0$.  To resolve it, one can use the stated 
%third order differential equation. A standard math package then solves it
%either numerically or expands it into a power series,
%%
%\bdm
%k(t)\ = \ t - \frac{1}{15}t^5 + \frac{13}{2835}t^9-\frac{1888}{6081075}t^{13}
%+ O(t^{15}).
%\edm
%%
Figure \ref{wind-geod} shows the shape of two typical geodesics of $\nabla$
for this vector field; the first goes through the origin, the second
through the point $(0,2)$. Observe that for the first geodesic, its
extension in negative time coincides with its development in positive time
rotated by $\pi$. The geodesic through  $(0,0)$ with zero slope resembles
an archimedean spiral, but turns out not to be one.

%---------------------------------------------------------------------------
\begin{figure}
\bdm
\begin{minipage}[b]{.47\linewidth}
\epsfig{figure=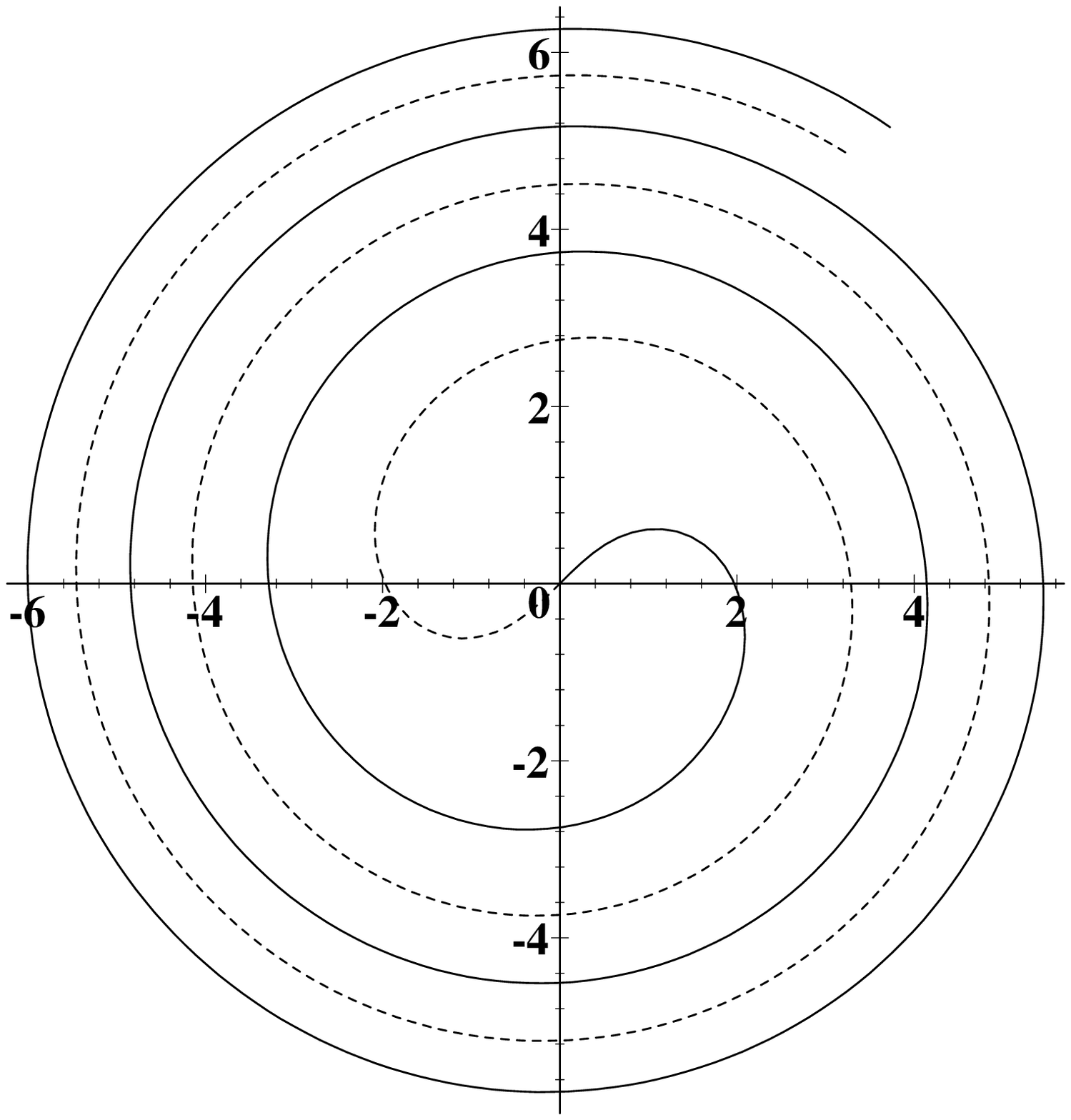,width=\linewidth, bbllx=50, bblly=120, bburx=580,bbury=690}
%\caption{geodesic through $(0,0)$ with slope $(1,1)$}
\end{minipage}\hfill
%---------------------------------
\begin{minipage}[b]{.47\linewidth}
%\centering
\epsfig{figure=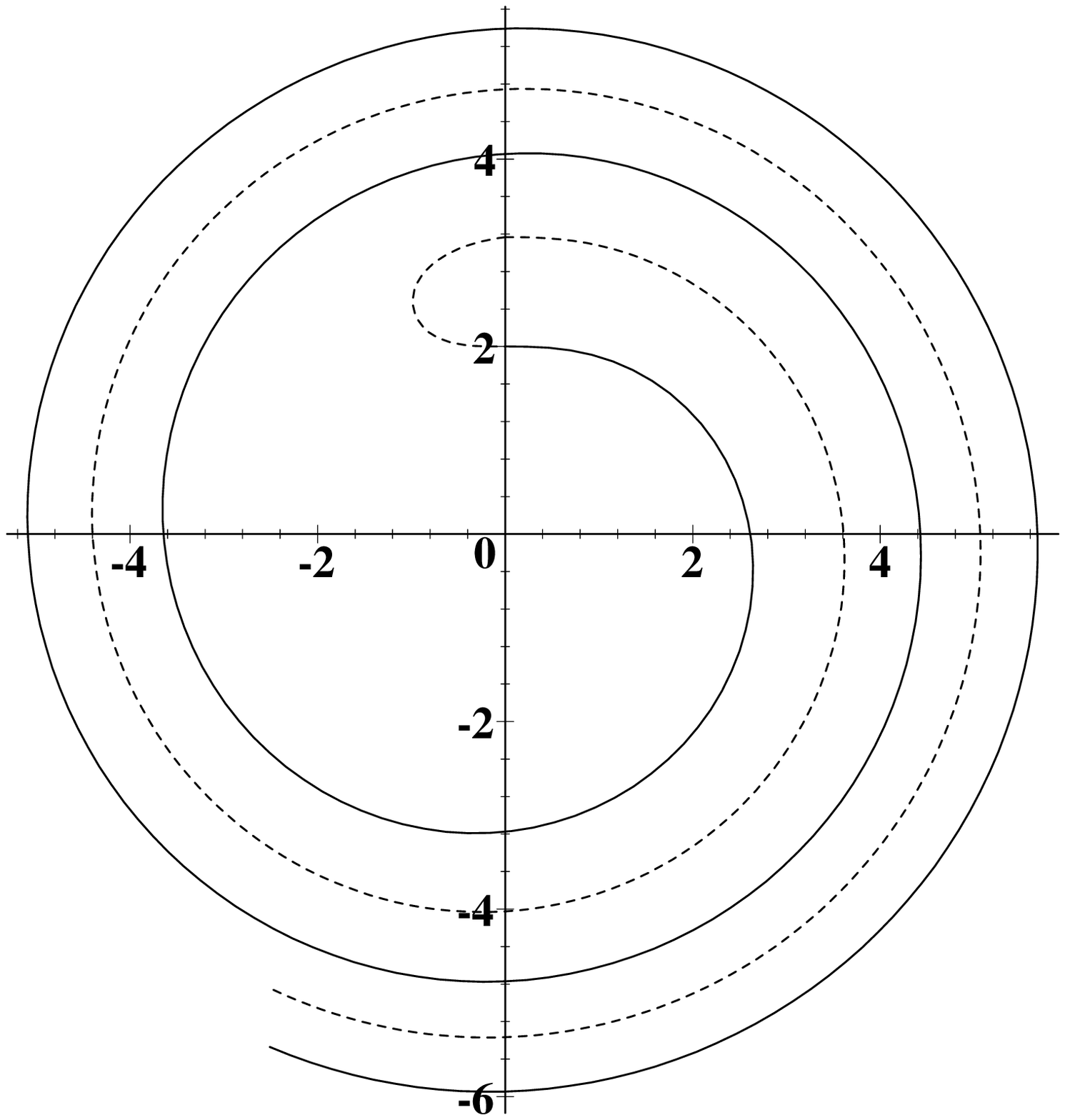, width=\linewidth, bbllx=50, bblly=140, bburx=580,bbury=690}
%\caption{geodesic through $(0,2)$ with slope $(1,0)$}
\end{minipage}
%\end{center}
\edm
\caption{Geodesic through $(0,0)$ with slope $(1,1)$ and through
$(0,2)$ with slope $(1,0)$ for $V=x\del_y-y\del_x$ ($t>0$ solid, 
$t<0$ dashed).}\label{wind-geod}
\end{figure}
%----------------------------------------------------------------------------
\end{exa}
\begin{exa}
%----------
The vector field $V=y\del_x$ ($f=y,g=0$) is obtained for $p=y^2/2$.
The geodesic equations read
\bdm
\ddot{x}+y\,\dot{y}^2\,=\,0,\quad \ddot{y}-y\,\dot{x}\dot{y}\,=\,0.
\edm
In this particular case, the second invariant of motion can be stated and
discussed rather explicitly.
For simplicity, set $E=1$. Inserting $\dot{x}=\pm \sqrt{1-\dot{y}^2}$
into the second geodesic equation, we obtain an ordinary differential
equation for $y$ alone, 
\bdm
\ddot{y}\ =\ \pm y \dot{y}\sqrt{1-\dot{y}^2}\ \text{ for }
\dot{x}\ =\ \pm \sqrt{1-\dot{y}^2}.
\edm
\begin{lem}
%----------
The ordinary differential equation $\ddot{y}=\pm y \dot{y}\sqrt{1-\dot{y}^2}$
has the invariant of motion $c:= \pm\, y^2/2 -\arcsin \dot{y}$. \qed
\end{lem}
\noindent
Formally, this allows us to fully integrate the differential equation,
\bdm
t+\tilde{c}\ =\ \int_{y_0}^y \frac{dy}{\sin(\pm y^2/2-c)}.
\edm
The integrand becomes singular whenever there exists a $k\in\Z$ such
that $\pm\,y^2=2c+2k\pi$. For a given $y_0$, two cases can be distinguished.
Either $c$ is such that $y_0$ is a zero of the denominator, then
the integral has to be interpreted in such a way that $t$ growths to
$\pm\infty$ for this value $y=y_0$. Thus, $y$ is the constant function $y_0$,
and one checks that indeed, $(at+b,y_0)$ is a solution of the geodesic
equations. Otherwise, a given $y_0$ lies in some interval between two
singularities $y_1,y_2$ of the integrand, $y_0\in (y_1,y_2)$; then,
$t$ will become infinitely large as $y$ approaches $y_1$ and $y_2$,
hence meaning that $y(\pm\infty)\ra y_1,y_2$. A generic geodesic lies thus
in a horizontal strip, which eventually can degenerate to a line.
Figure \ref{Hopf-Rinow-counterexa} illustrates two typical geodesics
of this connection. The vector field $V$ commutes with $\del_x$, hence
horizontal translates of geodesics are again geodesics by 
Lemma~\ref{isometries}.
\begin{NB}
%---------
Hence, the Hopf-Rinow theorem does not necessarily hold for
metric connections with vectorial torsion: for $V=y\del_x$, two
points lying in different ``geodesic strips'' cannot be joined by a geodesic 
arc. 
\end{NB}
%
%----------------------------------------------------------------------------
\begin{figure}
\bdm
\begin{minipage}[b]{.47\linewidth}
\epsfig{figure=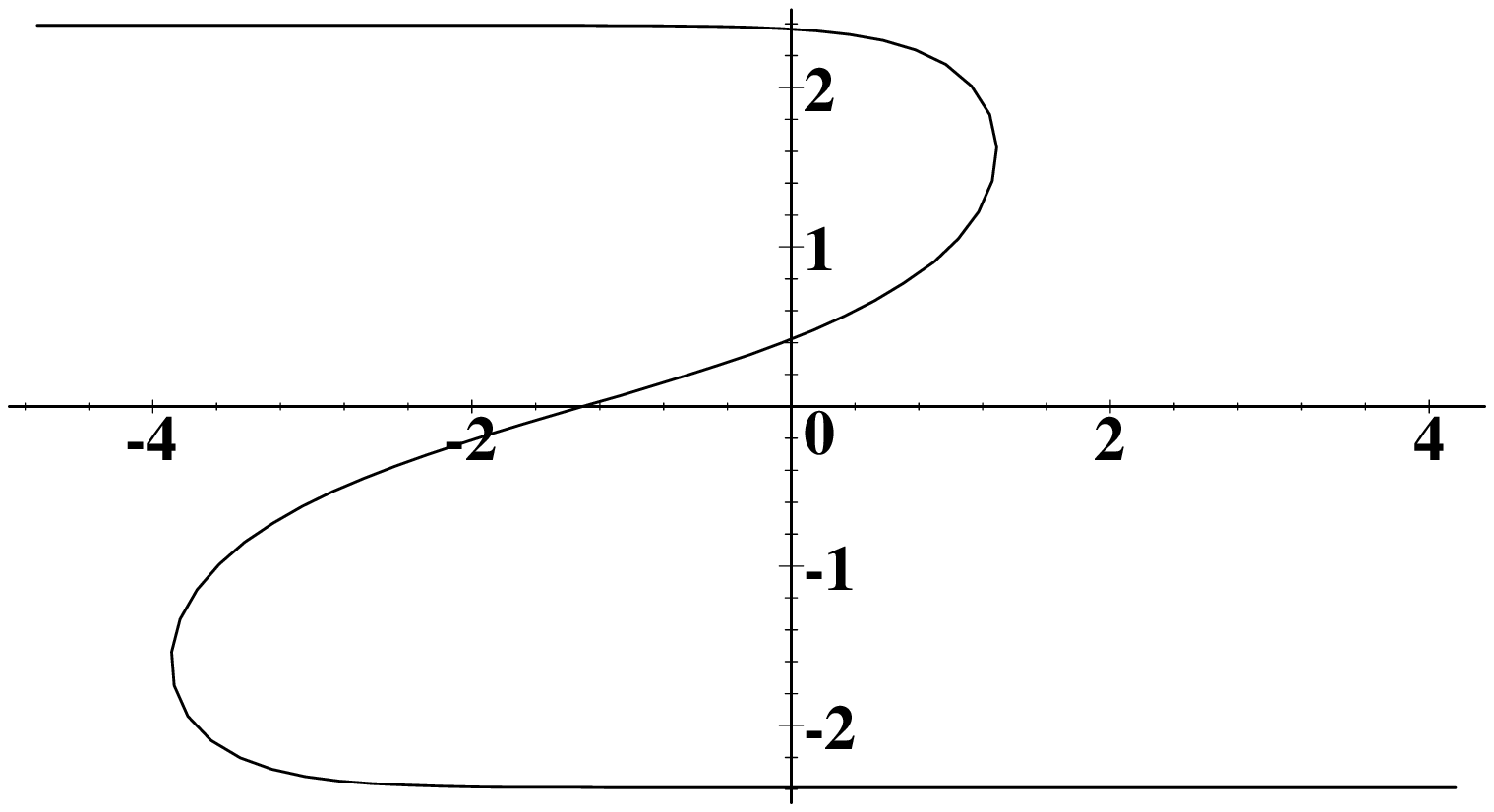,width=\linewidth, bbllx=50, bblly=230, bburx=560,bbury=590}
%\caption{geodesic through $(0,0)$ with slope $(1,1)$}
\end{minipage}\hfill
%---------------------------------
\begin{minipage}[b]{.47\linewidth}
%\centering
\epsfig{figure=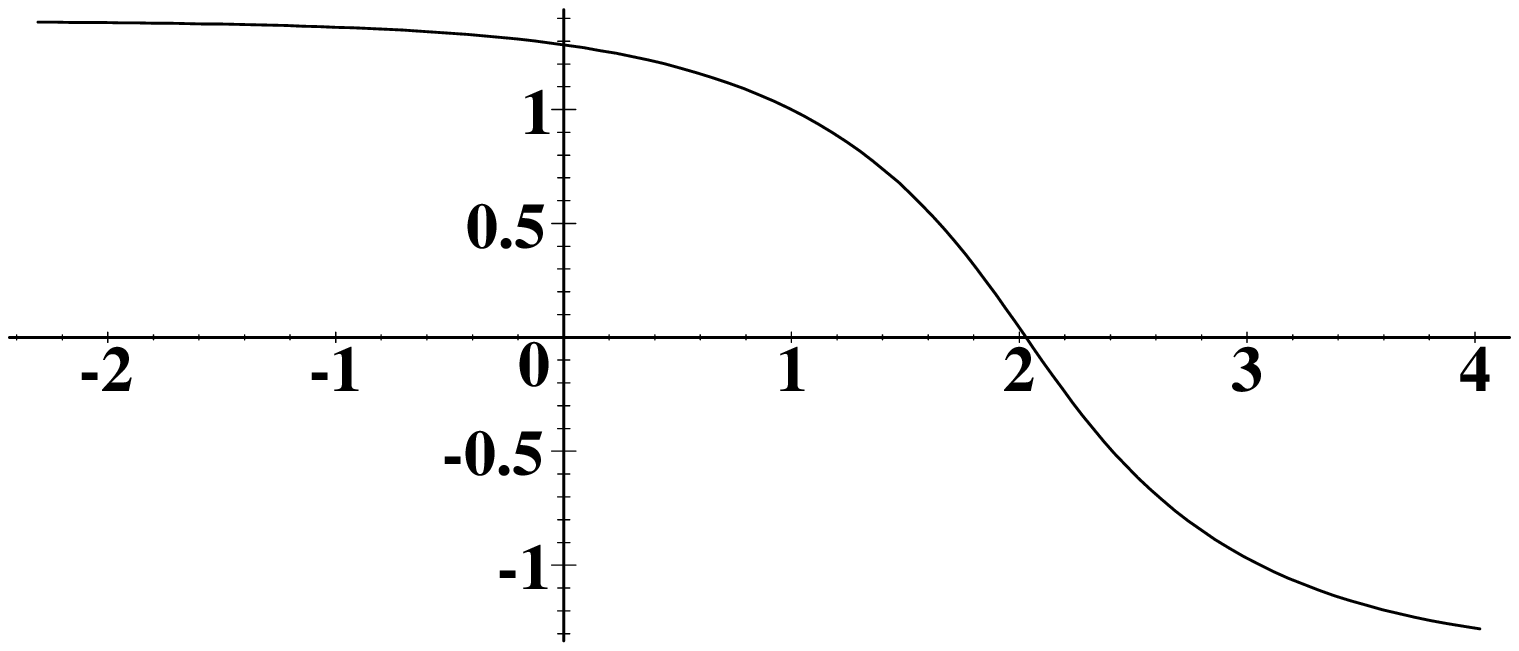, width=\linewidth, bbllx=50, bblly=230, bburx=560,bbury=590}
%\caption{geodesic through $(0,2)$ with slope $(1,0)$}
\end{minipage}
%\end{center}
\edm
\caption{Geodesics through $(1,1)$ with slopes $(1,1)$ and $(-1,1/2)$ for 
the vector field $V=y\del_x$.}\label{Hopf-Rinow-counterexa}
\end{figure}
%------------------------------------------------------------------------
\end{exa}

%--------------------------------------------------------------------------
%\nocite{Feeman}\nocite{Strubecker2}
%
%\bibliographystyle{amsalpha}
%---------------------------
%\bibliography{/home/agricola/bib/preamble,/home/agricola/bib/ia-pub,/home/agricola/bib/mathbooks,/home/agricola/bib/mathpapers}       
%--------------------------------------------------------------------------
\providecommand{\bysame}{\leavevmode\hbox to3em{\hrulefill}\thinspace}

%----------------------------------------------------------------------------
\end{document}